\documentclass[12pt]{amsart} 
\usepackage{amsmath,amscd,amssymb,graphics,color,epsf}
\usepackage[all]{xy}

\newcommand{\Z}{\mathbb{Z}}
\newcommand{\R}{\mathbb{R}}
\newcommand{\Q}{\mathbb{Q}}

\newcommand{\PP}{\mathbb{P}}

\newcommand{\F}{\mathcal F} 
\newcommand{\Sym}{\operatorname{Sym}}

\newcommand{\<}{\langle}   
\renewcommand{\>}{\rangle} 

\begin{document}

\title{Tropical Abelian varieties, Weil classes and the Hodge Conjecture}
\author{Ilia Zharkov}
\address{Kansas State University, 138 Cardwell Hall, Manhattan, KS 66506}
\email{zharkov@ksu.edu}

\begin{abstract}
We describe in some details an idea of M.~Kontsevich how one can try to find a counterexample to the Hodge conjecture using tropical geometry.
\end{abstract}

\maketitle

Let $\Gamma_1, \Gamma_2$ be two lattices in $V\cong \mathbb R^g$ and $Q:\Gamma_2^*\to \Gamma_1$ be a lattice isomorphism, called the {\em polarization}. We say that $X=V/\Gamma_1$ is a {\em principally polarized tropical abelian variety} if the polarization $Q$ is symmetric and positive definite, viewed as a bilinear form on $V$. Any such tropical abelian variety $X$ arises as a maximal degeneration of a family of complex abelian varieties $\mathcal A$ where half of the 1-cycles vanish, they form a lattice $\Gamma_2$, and the remaining ones form the lattice $\Gamma_1$. 

The tropical homology cosheaf $\F_p$, see \cite{IKMZ}, in this case is constant and it can be canonically identified with $\wedge^p \Gamma_2$ (the $p$-th homology of the collapsed torus fibers). Then the tropical homology $H_q (X, \F_p)$ is identified with $\wedge^q \Gamma_1 \otimes \wedge^p \Gamma_2$.

Among elements in $H_p(X, \F_p)$ we can single out the tropical Hodge $(p,p)$-classes as the kernel of the eigenwave action map, see \cite{MZ13}:
$$\phi: \wedge^p \Gamma_1 \otimes \wedge ^p \Gamma_2 
\to \wedge^{p-1} \Gamma_1 \otimes \wedge ^{p+1} V.$$
On the other hand, the algebraic tropical cycles are (weighted) balances polyhedral complexes with $\Gamma_2$-rational slopes. To every such algebraic cycle $Z$ one can associate the {\em tautological} $(p,p)$-chain $vol (Z)$ by taking the support of $Z$ and framing every (oriented) cell in it with its integral volume element in $\wedge^p \Gamma_2$. The balancing condition guaranties that this chain is a polyhedral $p$-cycle with coefficients in $\wedge ^p \Gamma_2 $, that is, it provides an element in $\wedge^p \Gamma_1 \otimes \wedge^p \Gamma_2$. Moreover a simple deformation argument shows that its class has to be in the kernel of $\phi$, i.e. $vol (Z)$ defines a Hodge class. The tropical Hodge conjecture says that every Hodge class is represented by an algebraic cycle.

If the tropical Hodge fails for such an $X=V/\Gamma_1$ then a generic member of the family of complex abelian varieties 
$$X_\epsilon=(\Gamma_2\otimes{\mathbb{C}}^*)/\epsilon^{-1}e^{\Gamma_1}, \ \epsilon\to 0,
$$
would provide a counterexample to the classical Hodge conjecture. The converse implication though may be false.

From now on we concentrate on dimension four, the lowest dimension to look for a counterexample to the Hodge conjecture. 
There are abelian varieties of so-called Weil type \cite{Weil} which have an additional symmetry of the period lattice (allowing ``complex multiplication'') resulting in two extra Hodge $(2,2)$-classes (called the Weil classes) in
$$(\wedge^2 \Gamma_1 \otimes \wedge^2 \Gamma_2 )\cap \ker \phi.$$
Classically the Weil classes are {\em not} known to represent algebraic cycles except in certain special cases. A good reference for the Hodge conjecture for abelian varieties is \cite{geemen}.

Following M.~Kontsevich \cite{kontsevich} we now give details on how one can explicitly disprove the tropical Hodge conjecture for abelian varieties of Weil type. Let us fix a basis $\{e_1, e_2, e_3, e_4\}$ of the lattice $\Gamma_2$, and also fix a positive integer $d$. Then the  4-dimension family of tropical Weil abelian varieties $X_{a,b,c,e}$ can be described by the polarization matrix  
$$
Q=\left( 
\begin{smallmatrix} 
a & b & 0 & e \\
b & c & -e & 0 \\
0 & -e & da & db \\
e & 0 & db & dc
\end{smallmatrix} 
\right), 
$$
where $a, b, c, e$ are real numbers such that $a>0$ and $d\,(ac-b^2)-e^2 > 0$.
The columns $(\gamma_1 \ \gamma_2 \ \gamma_3 \ \gamma_4)$ of the matrix $Q$ can be thought of as the coordinates of a basis for the other (period) lattice  $\Gamma_1\subset \Gamma_2\otimes \R$. We denote the parameter lattice $\Z\<a,b,c,e\> \cong \Z^4$ by $\Gamma_p$.
Then $\Gamma_1$ can be thought of as a rank 4 sublattice in the rank 16 lattice $\Gamma_2\otimes \Gamma_p$.

The ``complex multiplication'' $(\sqrt{-d})_*$  acts as follows:
\begin{gather*}
(e_1, e_2, e_3, e_4) \mapsto (de_3, de_4, -e_1, -e_2)\\
(\gamma_1 , \gamma_2 , \gamma_3 , \gamma_4) \mapsto 
  (\gamma_3 , \gamma_4 , -d \gamma_1 , -d \gamma_2).
\end{gather*}

We use the notations $\gamma_{ij}=\gamma_i \wedge \gamma_j$ and $e_{ij}=e_i \wedge e_j$. 
The square of the polarization
$$ \theta = \sum_{i < j} \gamma_{ij} \otimes e_{ij} \in \wedge^2 \Gamma_1 \otimes \wedge^2 \Gamma_2
$$
 is always a Hodge (2,2)-class. The Weil classes $w_1, w_2 \in H_2(X,\mathcal F_2) = \wedge^2 \Gamma_1 \otimes \wedge^2 \Gamma_2$ are:
\begin{multline*}
w_1=
\gamma_{12} \otimes e_{12} 
-\frac1d \gamma_{34} \otimes e_{12} 
-\gamma_{14} \otimes e_{14} 
-\gamma_{14} \otimes e_{32} \\
-\gamma_{32} \otimes e_{14} 
-\gamma_{32} \otimes e_{32}
-d \gamma_{12} \otimes e_{34} 
+\gamma_{34} \otimes e_{34}, 
\end{multline*}
and
\begin{multline*}
w_2=
\gamma_{14} \otimes e_{12}
- d \gamma_{14} \otimes e_{34} 
+ d \gamma_{12} \otimes e_{32} 
- \gamma_{34} \otimes e_{32} \\
+ d \gamma_{12} \otimes e_{14} 
- \gamma_{34} \otimes e_{14} 
+ \gamma_{32} \otimes e_{12} 
- d \gamma_{32} \otimes e_{34}. 
\end{multline*}
By considering $\Gamma_1$ as a sublattice in $\Gamma_2\otimes \Gamma_p$ we can write the above classes  in $\Sym^2 \Gamma_p \otimes \Sym^2 (\wedge^2 \Gamma_2)$ (a straightforward calculation) as:

\begin{multline*}
\theta = d ( -\frac{e}{d} e_{12} +a e_{13}+ b (e_{14}+ e_{23}) +c e_{24} - e e_{34}) ^2 \\
 + D ( -\frac 1d e_{12}^2 + e_{14}^2 + e_{23}^2 - 2 e_{13} e_{14} +d e_{34}^2),
\end{multline*}
and
\begin{align*}
w_1 & = D (\frac1{\sqrt{d}} \, e_{12}+\sqrt{d} \, e_{34})^2 - D (e_{14}+e_{32})^2,\\
w_2 & = 2 D (e_{12}-d \, e_{34})(e_{14}+e_{32}),
\end{align*}
where $D= d\,(ac-b^2)-e^2 \in \Sym^2 \Gamma_p$.

M. Kontsevich \cite{kontsevich} suggested the following scheme of showing that, for a generic choice of parameters $a,b,c,e$, the Weil classes are {\em not} represented by tropical algebraic cycles. First, observe that a {\em family} of tropical cycles in the family $X_{a,b,c,e}$ will have vertices which vary rationally over the space of parameters $(a, b, c, d)$. That is we can assume that our tropical cycles have vertices in $\Gamma_2 \otimes \Gamma_p$, which makes the setup countable. 

Next we consider two spaces. The first consists of all polygonal chains in $X_{a,b,c,e}$ with vertices in $\Gamma_2 \otimes \Gamma_p /\Gamma_1$ and edges having $\Gamma_2$-rational slopes. Any polygon with rational slope edges may be subdivided into triangles and parallelograms whose edges also have rational slopes. Thus we can assume that all chains consist only of triangles and parallelograms.
One subtlety here is that 
we may need to tensor everything with $\Q$. The reason is when we subdivide a polygon with vertices in the lattice $\Gamma_2 \otimes \Gamma_p$ into triangles and parallelograms the new vertices may be only rational. Another way to resolve this problem is to refine the lattice by allowing at every computational step a finite set of denominators (which still leaves the setup countable).

The second space is the space of affine flags (vertex $\subset$ edge $\subset$ 2-face) $(F_0,  F_1, F_2)$ where 
 $F_0$ is a vertex in $(\Gamma_2 \otimes \Gamma_p )/\Gamma_1$, $F_1 \in \PP(\Gamma_2\otimes \Q)$  is the direction of the edge and  $F_2 \in \wedge^2 \Gamma_2$ is the unit volume element of the 2-face. We add the flags formally and only use the linear structure on the $F_2$-component provided that the flags have common first two $F_0, F_1$.
  
Consider the map $\alpha$ from the space of polygonal chains to the flags which sends each polygon to the sum of all its flags. Note that the balancing condition says that any tropical cycle is mapped by $\alpha$ to the combination of flags whose $F_2$-components are all 0. On the other hand, there is the tautological map from the polygonal chains to $\Sym^2 \Gamma_p \otimes \Sym^2 (\wedge^2 \Gamma_2) $ which sends each polygon to its area form equipped with the unit area element (see the RHS's of \eqref{eq:triangle} and \eqref{eq:parallelogram} below for triange and parallelogram). This map coincides with the tautological map $Z\mapsto vol(Z)\in \wedge^p \Gamma_1 \otimes \wedge^p \Gamma_2$ on tropical cycles (hence the same name) if one identifies $\Gamma_1$ as a sublattice in $\Gamma_2\otimes \Gamma_p$.

The idea of Kontsevich is to build a map $\Phi$, linear on $\wedge^2 \Gamma_2\otimes \R$,  such that the diagram
 $$\xymatrix{
 \left\{ \txt{chains of $\triangle$'s and $\lozenge$'s\\
 with  $\Gamma_2$-integral slopes\\ 
 with vertices in\\  
 $\Gamma_2 \otimes \Gamma_p$
 }\right\} 
 \ar[d]_--{vol} \ar[rr]^-\alpha & &
   \left\{  \txt{ 
 Affine flags $(F_0, F_1, F_2)$ \\
  $F_0\in (\Gamma_2 \otimes \Gamma_p )/\Gamma_1$\\
  $F_1\in \PP(\Gamma_2\otimes \Q)$ \\ $F_2\in \wedge^2 \Gamma_2$} 
  \right\}
  \ar[lld]^{\Phi} \\
 \Sym^2 \Gamma_p \otimes \Sym^2 (\wedge^2 \Gamma_2) 
 }
 $$
commutes modulo a {\em proper} sublattice of $\Z\< \theta, w_1, w_2\>$ in $ \Sym^2 \Gamma_p \otimes \Sym^2 (\wedge^2 \Gamma_2) $. Then there will be a Hodge class $w\in \Sym^2 \Gamma_p \otimes \Sym^2 (\wedge^2 \Gamma_2) $ which is {\em not} represented by any algebraic cycle $Z$ under the map $vol$ since all such cycles are killed by the composition $\Phi \circ \alpha$.


Here is an explicit setup for computer programming. We need a family of linear maps 
$$\Phi_{x,u}: \wedge^2 \Gamma_2 \to \Sym^2 \Gamma_p \otimes \Sym^2 (\wedge^2 \Gamma_2)
$$ 
parameterized by two elements: 
$x\in (\Gamma_2 \otimes \Gamma_p )/\Gamma_1 \cong \Z^{12}$ and $u\in \PP(\Gamma_2\otimes \Q)\cong \Q\PP^3$ and which satisfies two sets of equations (commutativity of the diagram). 

The first type (triangles): for every $x\in (\Gamma_2 \otimes \Gamma_p )/\Gamma_1$,  $s\in \Gamma_p$ and $u,v \in \Gamma_2$:
\begin{multline}\label{eq:triangle}
(\Phi_{x,u}-\Phi_{x,v}+\Phi_{x+su,u-v}-\Phi_{x+su,u}+\Phi_{x+sv,v}-\Phi_{x+sv,u-v})(u\wedge v) \\
= s^2\otimes (u\wedge v)^2. 
\end{multline}
Here and later we will abuse the notations for $u,v \in \Gamma_2$ to also denote their images in 
$\PP(\Gamma_2\otimes \Q)$ in the subscript of $\Phi$.

The second type (parallelograms): for every $x\in (\Gamma_2 \otimes \Gamma_p )/\Gamma_1$,  $s, t \in \Gamma_p$ and $u,v \in \Gamma_2$:
\begin{multline}\label{eq:parallelogram}
(\Phi_{x,u}-\Phi_{x,v}+\Phi_{x+su,v}-\Phi_{x+su,u} - \Phi_{x+tv,u} + \Phi_{x+tv,v}\\
 +\Phi_{x+su+tv,u} - \Phi_{x+su+tv,v})(u\wedge v) = 2 st (u\wedge v)^2.
\end{multline}

Since the class $\theta$ is algebraic, we know that it is impossible to solve these equations on the nose. Our goal is to solve them modulo some proper sublattice of $\Z\< \theta, w_1, w_2\>$. 

The following ansatz solves the equations modulo $\Z\< \theta, w_1, w_2\>$. We assume the differences
$$\lambda_{x,s,u}= \Phi_{x+su, u} - \Phi_{x,u} 
$$ 
to be linear in $s\in \Gamma_p$ and in $u\in \Gamma_2$. Together with the obvious cocycle condition this forces 
$$\lambda_{x+tu, s, u} = \lambda_{x,s,u},
$$
leaving only finitely many degrees of freedom. Then the problem is to make sure that the $\lambda_{x, \cdot, \cdot}$ are defined modulo $\Gamma_1$. This gives rise to a finite linear system which can be solved explicitly, but unfortunately the solution does not hold modulo any proper sublattice of $\Z\< \theta, w_1, w_2\>$.

\
\bibliographystyle{plain}

\end{document}